\documentclass[11pt]{article}
\usepackage{graphicx} 
\usepackage{fullpage}
\usepackage{amsfonts}
\usepackage{graphicx}
\usepackage{float}
\usepackage{tikz}
\usepackage{setspace} 
\usepackage{xcolor}

\usepackage[pdfdisplaydoctitle=true]{hyperref}

\newcommand{\inv}{\mathrm{inv}}

\usepackage{amssymb}
\usepackage{epstopdf}
\usepackage{framed}
\usepackage{fullpage}

\usepackage{setspace}
\usepackage{amsmath,amsthm}
\newtheorem{thm}{Theorem}
\newtheorem{proposition}{Proposition}

\newtheorem{definition}{Definition}

\theoremstyle{remark}
\newtheorem{remark}{Remark}

\title{Asymptotics of dynamic ASEP using duality}
\author{Jeffrey Kuan and Zhengye Zhou }
\date{}
\begin{document}

\Large 

\maketitle

\begin{center}
\fbox{\parbox{0.8\textwidth}{
Accessibility Statement: A WCAG2.1AA compliant version of this PDF will eventually be available at the first author's professional webpage \cite{TSB1}. 

\vspace{0.1in}

Some guidelines of the British Dyslexia Association Style Guide were consulted in creating this document. The dyslexia--friendly color contrast is on the external link.

}
}
\end{center}

\abstract{Using a recently developed method for proving asymptotics via orthogonal polynomial duality \cite{KZ23}, we prove that the dynamic ASEP introduced in \cite{Bor20} has asymptotics which are either distributed as the Tracy--Widom \(F_2,\) or are almost surely bounded. Using a different duality, we also provide contour integrals formulas for multi--species ASEP, which generalize results for the single--species ASEP.}

\tableofcontents

\Large

\section{Introduction}

In recent years, the authors have developed a method \cite{KZ23} to prove asymptotic fluctuations of models which have so--called ``orthogonal polynomial duality.'' (See \cite{CFGGR,ACR21,CFG21,FG19,Gro19,REU2020,ZZ21,FKZ22,FRS22,REU2022,GW23} for papers on the topic of orthogonal polynomial duality). The essential idea is that a desired observable, such as the height function, decomposes over an orthogonal basis consisting of duality functions. This basis is indexed by the state space of the dual process; and in the context of interacting particle systems, are orthogonal with respect to certain ``blocking'' measures. One can then analyze the height function through this decomposition.

This paper will apply that method to so--called ``dynamic'' models, introduced in \cite{Bor20} and studied further in \cite{Agg18,BC20,CGM20}. In these dynamic asymmetric models, the asymmetry parameter ``dynamically'' changes depending on the height function. In the case of the dynamic ASEP, the asymmetry reverses as the particles drift. Although it is well established by now that for models in the KPZ universality class with step initial conditions, the asymptotic fluctuations will be Tracy--Widom \(F_2,\) there had not even been conjectures for the asymptotic fluctuations of asymmetric dynamic models (see \cite{Agg18} and \cite{BC20} for asymptotics of symmetric dynamic models). In \cite{KZ23}, the authors prove that the dynamic stochastic six vertex model with step initial conditions has Tracy--Widom \(F_2\) fluctuations, just as in the non--dynamic stochastic six vertex model.

In this paper, we prove that the dynamic ASEP with step initial conditions either has \(F_2\) fluctuations or are almost surely bounded, depending on the value of the dynamic parameter. Based on these results, it would seem that dynamic asymmetric models have the same asymptotic fluctuations as the usual asymmetric models with asymmetry parameter \(q\) or \(q^{-1}.\) There is a subtlety here to be noted: for ASEP, inverting the asymmetry parameter reverses the direction of the drift; in contrast, the stochastic six vertex model is totally asymmetric, and inverting the asymmetry parameter does not reverse the direction of the drift.

Our method will utilize duality functions found by \cite{GW23}, which discovered duality functions between dynamic (generalized) ASEP and a usual (generalized) ASEP. The duality functions are written in terms of \(q\)--Hahn orthogonal polynomials, which degenerate to duality functions expressed in terms of quantum \(q\)--Krawtchouk orthogonal polynomials. The latter duality functions are duality functions for the usual (generalized) ASEP. Using the previously found methods of \cite{KZ23} and estimates on the duality function and blocking measures, this allows one to compare the asymptotics of the dynamic ASEP to the usual ASEP. 

Additionally, we will present some contour integral formulas for certain observables in multi--species ASEP (introduced in \cite{Ligg76}). These formulas are derived using contour integral formulas in \cite{KuanAHP} and dualities of \cite{BS3,Kuan-IMRN} (see also \cite{Kua16,BS,BS2} for the two--species case; and \cite{Sch97} for the one--species case). These formulas generalize formulas found those in \cite{TW08}.

\section{Background}

\subsection{\(q\)--notation and orthogonal polynomials}

In this subsection, we define \(q\)--notation and some orthogonal polynomials. 

Recall the $q$--Pochhammer 
$$(a ; q)_{n}=(1-a)(1-a q) \cdots\left(1-a q^{n-1}\right), \quad n \in \mathbb{N}
$$
and the $q$--hypergeometric series
$$
{ }_{r+1} \varphi_{r}\left(\begin{array}{c}a_{1}, \ldots, a_{r+1} \\ b_{1}, \ldots, b_{r}\end{array} ; q, z\right)=\sum_{n=0}^{\infty} \frac{\left(a_{1} ; q\right)_{n} \cdots\left(a_{r+1} ; q\right)_{n}}{\left(b_{1} ; q\right)_{n} \cdots\left(b_{r} ; q\right)_{n}} \frac{z^{n}}{(q ; q)_{n}}.
$$
If some $a_k=q^{-l}$ for some non--negative integer $l$, then the series will terminate because $(q^{-l};q)_{l+1}=0$. It will sometimes be helpful to define the \( q\)--binomial is:
\begin{equation}
\left[\begin{array}{l}n \\ k\end{array}\right]_{q}=\frac{(q ; q)_{n}}{(q ; q)_{k}(q ; q)_{n-k}}=\frac{\left(q^{-n} ; q\right)_{k}}{(q ; q)_{k}}\left(-q^{n}\right)^{k} q^{-\frac{1}{2} k(k-1)}.
\end{equation}

The \(q\)--Hahn polynomial (section 14.6 of \cite{KLS}) is defined by
\[
P_{n}\left(q^{-x} ; \alpha, \beta, N \mid q\right)= \ _{3} \varphi_{2}\left(\begin{array}{c}q^{-n}, \alpha \beta q^{n+1}, q^{-x} \\ \alpha q, q^{-N}\end{array} ; q, q\right), \quad n=0,1,2, \ldots, N
\]
and satisfy the orthogonality relations
\Large
\[
\begin{array}{l}\sum_{x=0}^{N} \frac{\left(\alpha q, q^{-N} ; q\right)_{x}}{\left(q, \beta^{-1} q^{-N} ; q\right)_{x}}(\alpha \beta q)^{-x} P_{m}\left(q^{-x} ; \alpha, \beta, N \mid q\right) P_{n}\left(q^{-x} ; \alpha, \beta, N \mid q\right) \\ =\frac{\left(\alpha \beta q^{2} ; q\right)_{N}}{(\beta q ; q)_{N}(\alpha q)^{N}} \frac{\left(q, \alpha \beta q^{N+2}, \beta q ; q\right)_{n}}{\left(\alpha q, \alpha \beta q, q^{-N} ; q\right)_{n}} \frac{(1-\alpha \beta q)(-\alpha q)^{n}}{\left(1-\alpha \beta q^{2 n+1}\right)} q^{\binom{n}{2}-N n} \delta_{m n}\end{array}
\]
The quantum \(q\)--Krawtchouk polynomial is defined by (section 14.14 of \cite{KLS})
\[
K_n^{\text{qtm}}(q^{-x};p,c;q) := \ _{2}\varphi_1(q^{-x},q^{-n};q^{-c};q,pq^{n+1})
\]
and is related to the \(q\)--Hahn polynomial by 
\[
\lim _{\alpha \rightarrow \infty} P_{n}\left(q^{-x} ; \alpha, p, N \mid q\right)=K_{n}^{q t m}\left(q^{-x} ; p, N ; q\right)
\]
and satisfy the orthogonality relations
\[
\begin{array}{l}\sum_{x=0}^{N} \frac{(p q ; q)_{N-x}}{(q ; q)_{x}(q ; q)_{N-x}}(-1)^{N-x} q^{\binom{x}{2}} K_{m}^{q t m}\left(q^{-x} ; p, N ; q\right) K_{n}^{q t m}\left(q^{-x} ; p, N ; q\right) \\ =\frac{(-1)^{n} p^{N}(q ; q)_{N-n}(q, p q ; q)_{n}}{(q, q ; q)_{N}} q^{\binom{N+1}{2}-\binom{n+1}{2}+N n} \delta_{m n}, \quad p>q^{-N}.\end{array}
\]
Also define the one--site duality function
\[ \begin{array}{c}k^{\mathrm{qtm}}(n, x ; \lambda, \rho, v, N ; q)=K_{x}^{\mathrm{qtm}}\left(n ; \hat{p}, N ; q^{2}\right), \\ \text{ where }\hat{p}=v^{-1} q^{\rho-\lambda-N-1}\end{array} \]

\subsection{Definitions of models}

First, we will define the models that will be discussed in this paper.

\subsubsection{Definition of dynamic ASEP}

Here, we will essentially copy the definition of dynamic ASEP used by \cite{GW23}. 

\begin{definition}
    
$\operatorname{ASEP}_{\mathrm{R}}(q,  \rho)$ is a continuous-time Markov jump process on the state space $X=\{0,1\}^{\mathbb{Z}}$ depending on parameters $q>0$ and $\rho \in \mathbb{R}$. Given a state $\xi=\left(\xi_k\right)_{k=-\infty}^\infty \in X$ we define the height function $\left(h_k^{+}\right)_{k=-\infty}^\infty$ by
$$
h_k^{+}= h_{k,\rho}^+(\xi)=\rho+k+\sum_{j=k}^\infty 2 \xi_j.
$$
The generator is given by
$$
L_{q, \rho}^{\mathrm{R}} f(\xi)=\sum_{k \in \mathbf{Z}} C_k^{\mathrm{R},+}(\xi)\left[f\left(\xi^{k, k+1}\right)-f(\xi)\right]+C_{k+1}^{\mathrm{R},-}(\xi)\left[f\left(\xi^{k+1, k}\right)-f(\xi)\right] .
$$
Then a particle on site $k$ jumps to site $k+1$ at rate
$$
C_k^{\mathrm{R},+}(\xi)=q^{-1} \frac{\left(1+q^{-2 h_{k+1}^{+}}\right)}{\left(1+q^{-2 h_{k+1}^{+}-2}\right)},
$$
and a particle on site $k$ jumps to site $k-1$ at rate
$$
C_k^{\mathrm{R},-}(\xi)=q \frac{\left(1+q^{ -2 h_{k}^{+}}\right)}{\left(1+q^{-2 h_k^{+}+2}\right)} .
$$
\end{definition} 

Note that we use the convention that the height function counts particles to the right, so therefore \textit{step initial conditions} will mean particles are initially located at the negative integers; otherwise the height function will be infinite.

The dynamic ASEP interpolates between an ASEP and a reversed ASEP. For example, if $\rho=\infty$ and $q>1$, then the (local) drift is to the left; while for $\rho=-\infty$ and $q>1$, the local drift is to the right. As particles drift to the right the height function increases, and the drift will move towards the leftward direction. Thus, the dynamic ASEP has a tendency to push the height function downwards as it increases. For the remainder of this paper, assume that $q>1$. Note that the \(q<1\) case can be obtained by symmetry; namely, replace \(\rho\) with \(-\rho\) and invert the lattice. 

\subsubsection{Definition of (multi--species) ASEP }
Although the main results of this paper only concern the single--species ASEP, the multi--species (also called the ``colored'') ASEP will be used as a tool in the proofs. In this case, it suffices to consider the ``rainbow'' case when there is at most one particle of each species. We will additionally use an uncommon notation, which was used in \cite{KuanAHP}.

In this case, the state space will consist of pairs $(\mathbf{x},\sigma)$, where 
$$
\mathbf{x} \in \mathcal{W}_N^+ = \{(x_1 > x_2 > \ldots > x_N): x_i \in \mathbf{Z}\} \subset \mathbf{Z}^N
$$
and $\sigma \in S_N$ is a permutation on \(N\) letters. In more familiar notation using occupation variables, we can define a map \(\iota\) by defining \(\eta :=\iota(\mathbf{x},\sigma)\) as
\[
\eta(x)
= 
\begin{cases}
 0, & \text{ if } x \notin \mathbf{x} \\
 \sigma(k), & \text{ if } x=x_k \in \mathbf{x}.
\end{cases}
\]
for any \(x \in \mathbb{Z}.\) Here, \( \eta(x)= j\) means that there is a particle of type (variously called species, color or class) \(j\) at lattice site \(x.\)

We now define two generators, \(L^{\pm}_{\text{rainbow}}.\) These are defined by having off--diagonal entries
\begin{multline*}
L^{\pm}_{\text{rainbow}} (( \mathbf{x},\sigma),(\mathbf{x}',\sigma')) \\
= 
\begin{cases}
1, \text{ if } \mathbf{x}' = \mathbf{x}^{x \rightarrow x+1} \text{ for some } x\in X \text{ and } \sigma=\sigma',\\
q, \text{ if } \mathbf{x}' = \mathbf{x}^{x \rightarrow x-1} \text{ for some } x\in X \text{ and } \sigma=\sigma',\\
1, \text{ if } \mathbf{x}' = \mathbf{x} \text{ and } \sigma' = \sigma \circ (r \ r+1) \text{ for some } r \text{ and } \inv(\sigma') = \inv(\sigma) \pm 1,\\ 
q, \text{ if } \mathbf{x}' = \mathbf{x} \text{ and } \sigma' = \sigma \circ (r \ r+1) \text{ for some } r \text{ and } \inv(\sigma') = \inv(\sigma) \mp 1,\\
0, \text{else}. 
\end{cases}
\end{multline*}
Here, the superscript \({x} \rightarrow x\pm 1\) in \(\mathbf{x}'\) indicates that \(\mathbf{x}'\) is obtained from \(\mathbf{x}\) by replacing \(\mathbf{x}\) with \(x \pm 1.\) The notation \( \mathrm{inv}(\sigma) \) indicates the number of inversions of \(\sigma,\) which is the number of pairs \((i,j)\) such that \(i<j\) and \(\sigma(i)>\sigma(j).\) The diagonal entries $L_{\text{rainbow}} (( \mathbf{x},\sigma),(\mathbf{x},\sigma)) $ are defined so that the rows sum to $0$. 

We describe this process in words. The set $\mathbf{x}$ indicates the locations where the sites are occupied by particles and $\sigma$ indicates the ordering of the species (or colors) particles. The right jump rates are \(1\) and the left jump rates are \(q.\) Particles with higher ``priority'' swap places with particles of lower ``priority'' by ignoring their existence. Such swaps change the number of inversions in \(\sigma\) by \(1\) or \(-1.\) The choice of \(\pm 1\) determines whether particles of type \(i\) have priority over particles of type \(j,\) or vice versa, for fixed \(i<j.\)

We also note that there is something called a ``color--blind'' Markov projection from colored (multi--species) models to the usual (single--species) model. More specifically, we have that for all \(\sigma \in S_N,\)
\[
 \sum_{\sigma' \in S_N} L^{\pm}_{\text{rainbow}} (( \mathbf{x},\sigma),(\mathbf{x}',\sigma')) = L_{\text{ASEP}}(\mathbf{x},\mathbf{x}')
\]
where \(L_{\text{ASEP}}\) is the generator of ASEP with right jump rates \(1\) and left jump rates \(q.\)

\subsection{Previously known results}

\subsubsection{Orthogonal Dualities for dynamic ASEP}
Proposition 4.1(i) of \cite{GW23} gives a duality between a dynamic ASEP and a usual ASEP, in terms of \(q\)--Hahn polynomials. In the usual ASEP, left jump rates are \(q\) and right jump rates are \(q^{-1}.\) Since we are assuming that \(q>1\) this means that the dual process has drift to the left. They have a more general result for a duality between two dynamic generalized ASEPs (where up to \(N_k\) particles may occupy the lattice site \(k\)), in terms of \(q\)--Racah polynomials, but for reasons of brevity we do not state that result here. 

Define the one--site duality function by (equation (4.2) of \cite{GW23})
\[
\begin{array}{c}p(n, x ; \lambda, \rho, v, N ; q)=c_{\mathrm{p}}(n, x ; \lambda, \rho, v, N ; q) P_{x}\left(n ; \alpha, \beta, N ; q^{2}\right), \\ (\alpha, \beta)=\left(-v q^{\rho+\lambda-N-1}, v^{-1} q^{\rho-\lambda-N-1}\right)\end{array}
\]
where the coefficient is found in appendix A:
\[
c_{\mathrm{p}}(n, x ; \lambda, \rho, v, N ; q)=v^{n} \frac{\left(-v q^{\rho+\lambda-N+1} ; q^{2}\right)_{x}\left(v q^{2 n-\rho+\lambda-N+1} ; q^{2}\right)_{N}}{q^{n(n+\rho+\lambda-N)}\left(v q^{-2 x-\rho+\lambda+N+1} ; q^{2}\right)_{x+n}}.
\]
The duality function on \(L\) sites is:
\[
P_{\mathrm{R}}^{v}(\eta, \xi)=\prod_{k=1}^{L} p\left(\eta_{k}, \xi_{k} ; h_{k-1,0}^{-}(\eta), h_{k+1}^{+}(\xi), v, N_{k} ; q\right)
\]
here \(\xi\) is the original dynamic model and \(\eta\) is the dual model.

The orthogonality relation is in Proposition 4.4 (of \cite{GW23}):
\[
\begin{array}{l}\sum_{\eta \in X} P_{\mathrm{R}}^{v}(\eta, \xi) P_{\mathrm{R}}^{v}\left(\eta, \xi^{\prime}\right) \omega^{\mathrm{p}}(|\eta|) w(\eta ; \vec{N} ; q)=\frac{\delta_{\xi, \xi^{\prime}}}{\omega_{\mathrm{R}}^{\mathrm{p}}(|\xi|) W_{\mathrm{R}}(\xi ; \vec{N}, \rho ; q)} \\ \sum_{\xi \in X} P_{\mathrm{R}}^{v}(\eta, \xi) P_{\mathrm{R}}^{v}\left(\eta^{\prime}, \xi\right) \omega_{\mathrm{R}}^{\mathrm{p}}(|\xi|) W_{\mathrm{R}}(\xi ; \vec{N}, \rho ; q)=\frac{\delta_{\eta, \eta^{\prime}}}{\omega^{\mathrm{p}}(|\eta|) w(\eta ; \vec{N} ; q)}\end{array}
\]
where the \(\omega\) weights are
\[
\begin{aligned} \omega^{p}(x) & =\frac{v^{-2 x} q^{x(2 x-1)}\left(-v q^{\rho-|\vec{N}|+1} ; q^{2}\right)_{x}}{\left(v q^{-\rho+2 x-|\vec{N}|+1} ; q^{2}\right)_{|\vec{N}|-x}} \\ \omega_{\mathrm{R}}^{p}(x) & =\frac{\left(v q^{-\rho-2 x+|\vec{N}|+1} ; q^{2}\right)_{x}}{\left(-v q^{\rho-|\vec{N}|+1} ; q^{2}\right)_{x}}\end{aligned}
\]
and the \(w,W\) weights are
\[
w(\eta,\vec{N};q) = q^{\sum_{k=1}^M \eta_kN_k} q^{-2\sum_{k=1}^M\sum_{j=1}^k \eta_kN_j} \prod_{k=1}^M  q^{\eta_k(\eta_k-N_k)} \left[ \begin{array}{c} N_k \\ \eta_k \end{array}\right]_{q^2},
\]
\[
W_{\mathrm{R}}(\xi ; q, \rho)  =\prod_{k\in \mathbb{Z}} W\left(\xi_k ; q, 1, h_{k+1}^{+}(\xi)\right).
\]
For this paper, we are considering only the dynamic ASEP, where $N=1$. Thus, the latter weights simplify to
$$
\begin{aligned}
W(x ; q, 1, \rho) & =\frac{1+q^{4 x+2 \rho-2 }}{1+q^{2 \rho-2}} \frac{\left(-q^{2 \rho-2 } ; q^2\right)_x}{\left(-q^{2 \rho+2} ; q^2\right)_x} \frac{q^{-x(2 \rho+1+x-2 )}}{\left(-q^{-2 \rho} ; q^2\right)_1}
\end{aligned}
$$

Next, we discuss the degeneration to the quantum \(q\)--Krawtchouk orthogonal polynomial. Define (note we take \(\rho\) to \(\infty\) because \(q>1\))
\[
K_{\mathrm{qtm}}^{v}(\eta, \xi)=\lim _{\rho \rightarrow \infty} \frac{\left(v^{-2} q^{-2 \rho}\right)^{|\eta|}}{c^{v}(|\eta|,|\xi| ; 0,0) C^{v}(|\eta|,|\xi| ; 0,2 \rho)} P_{\mathrm{R}}^{v q^{\rho}}(\eta, \xi)
\]
then 
\[
K_{\mathrm{qtm}}^{v}(\eta, \xi)=\prod_{k=1}^{L} k^{\mathrm{qtm}}\left(\eta_{k}, \xi_{k} ; h_{k-1,0}^{-}(\eta), h_{k+1,0}^{+}(\xi), v, N_{k} ; q\right)
\]
where
\begin{multline*}
c(\zeta, \xi):=c^{v}(|\zeta|,|\xi| ; \lambda, \rho)\\
=\prod_{k=1}^{L} \frac{\left(v q^{2 \zeta_{k}-h_{k+1}^{+}(\xi)+h_{k-1}^{-}(\zeta)-N_{k}+1} ; q^{2}\right)_{N_{k}}}{\left(v q^{-2 \xi_{k}-h_{k+1}^{+}(\xi)+h_{k-1}^{-}(\zeta)+N_{k}+1} ; q^{2}\right)_{\xi_{k}+\zeta_{k}}}=\frac{\left(v q^{\lambda-\rho+2|\zeta|-|\vec{N}|+1} ; q^{2}\right)_{|\vec{N}|-|\zeta|}}{\left(v q^{\lambda-\rho-2|\xi|+|\vec{N}|+1} ; q^{2}\right)_{|\xi|}}
\end{multline*}
and
\[
C^{v}(|\zeta|,|\xi| ; \lambda, \rho)=\prod_{k=1}^{L} \frac{\left(-v q^{h_{k+1}^{+}(\xi)+h_{k-1}^{-}(\zeta)-N_{k}+1} ; q^{2}\right)_{\xi_{k}}}{\left(-v q^{h_{k+1}^{+}(\xi)+h_{k-1}^{-}(\zeta)-N_{k}+1} ; q^{2}\right)_{\zeta_{k}}}=\frac{\left(-v q^{\lambda+\rho-|\vec{N}|+1} ; q^{2}\right)_{|\xi|}}{\left(-v q^{\lambda+\rho-|\vec{N}|+1} ; q^{2}\right)_{|\zeta|}}
\]
The degenerations of the weights are:
\[
\begin{array}{l}\omega^{\mathrm{qtm}}(x)=v^{-x} q^{x(x+|\vec{N}|-1)}\left(v q^{2 x-N+1} ; q^{2}\right)_{|\vec{N}|-x} \\ \omega_{\mathrm{R}}^{\mathrm{qtm}}(x)=\frac{v^{x} q^{x(|\vec{N}|-x+1)}}{\left(v q^{1+|\vec{N}|-2 x} ; q^{2}\right)_{x}}\end{array}
\]

\subsubsection{Dualities and particle positions for multi--species ASEP}

A previous result of \cite{KuanAHP} calculates formulas for so--called ``$q$--exchangeable'' particle distributions in multi-species ASEP in terms of combinatorics of the symmetric group and contour integral formulas. 

Given $q$--exchangeable initial conditions supported at $\mathbf{y}$, the ``rainbow'' ASEP on $\mathbf{Z}$ satisfies
\begin{multline*}
\mathrm{Prob}((X,\sigma);t) \\
= \left( \frac{1}{2\pi i}\right)^N \frac{q^{\inv(\sigma)}}{N]_q^!} \sum_{\sigma \in S_N} \int_{\mathcal{C}_r} \cdots \int_{\mathcal{C}_r} A_{\sigma} \prod_{i} \xi_{\sigma(i)}^{x_i - y_{\sigma(i)}-1} e^{(1+q)\sum_i \epsilon(\xi_i)t} d\xi_1 \cdots d\xi_N,
\end{multline*}
where \( \epsilon(\xi)=\xi^{-1} + q\xi -1 \) and
$$
A_{\sigma} = \prod_{(\beta,\alpha) \text{ is an inversion of } \sigma} S_{(\beta,\alpha)}.
$$
with
$$
S_{\alpha,\beta}= - \frac{1+q\xi_{\alpha}\xi_{\beta} - \xi_{\alpha}}{1+q\xi_{\alpha}\xi_{\beta} - \xi_{\beta}}.
$$

There is also a self--duality function for the multi--species ASEP, found in [Kuan, Thm 2.5(b) and Proposition 5.2]; see also and [Belitsky-Schutz]. We write the dual process in terms of occupation variables, so the finite sets \(A^1,\ldots,A^n \subseteq \mathbb{Z}\) denote the locations of the species \(1,2,\ldots, n\) particles. This duality is ``triangular'' in the sense that it can be expressed as a triangular matrix for some indexing of the state spaces. First, let us define the set \(I\) on which the duality function is nonzero. The set \(I\) will consists of pairs of states of multi--species ASEP. Let \(\mathbf{A} = (A^1,\ldots,A^n)\) denote a state in the dual process. Then the set \(I\) is defined by setting
\[
I = \{ \left((\mathbf{x},\sigma), \mathbf{A} \right): \forall j \in [1,n] \text{ and } \forall x \in A^j, x=x_k \text{ for some } k \text{ and } \sigma(k) \geq j\}.
\]
In words, this means that for every species \(j\) particle in \( \mathbf{A},\) the corresponding lattice site in \( (\mathbf{x},\sigma)\) has a particle of species \( \geq j.\) The duality function is then
\[
\widehat{D}((\mathbf{x},\sigma),\mathbf{A}) = 
1_{ ((\mathbf{x},\sigma),\mathbf{A})\in I } \prod_{j=1}^n \prod_{z \in A^j} q^{-2z + 2N^{j,-}_z( (\mathbf{x},\sigma))  }
\]
where \( N_z^{j,-}((\mathbf{x},\sigma))\) is the number of particles of species \(\geq j\) to the left of \(z.\) In symbols,
\[
N_z^{j,-}((\mathbf{x},\sigma)) = \left| \{y \leq z: y =y_k \in \mathbf{x} \text{ where } \sigma(k) \geq j \}\right|.
\]

\subsubsection{Tracy--Widom contour integral formulas}

Theorem 3.1 of \cite{TW08} gives a contour integral formula for the left--most particle of ASEP when there are exactly \(N\) particles. That paper uses \(\alpha\) for the left jump rates and \(\beta\) for the right jump rates, normalized to \(\alpha+\beta=1.\)  To match with the notation here, let $\alpha=q/(q^{-1}+q)$ and $\beta = q^{-1}/(q^{-1}+q)$; note that $\alpha/\beta =q ^2$.

Assume initial conditions \(Y=(y_1<\cdots < y_N).\) Then the distribution of the left--most particle is
$$
\mathbb{P}\left(x_{1}(t)=x\right)=q^{N(N-1) / 2} \int_{\mathcal{C}_{r}} \cdots \int_{\mathcal{C}_{r}} I(x, Y, \xi) d \xi_{1} \cdots d \xi_{N}
$$
where
$$
I(x, Y, \xi)=\prod_{i<j} \frac{\xi_{j}-\xi_{i}}{q+q^{-1} \xi_{i} \xi_{j}-\xi_{i}} \frac{1-\xi_{1} \cdots \xi_{N}}{\left(1-\xi_{1}\right) \cdots\left(1-\xi_{N}\right)} \prod_{i}\left(\xi_{i}^{x-y_{i}-1} e^{\varepsilon\left(\xi_{i}\right) (1+q)t}\right)
$$
and $\varepsilon(\xi)=q \xi^{-1}+q^{-1} \xi-1$ and $\mathcal{C}_r$ are circles centered at the origin with radius $r<1$.

In equation (1.6) they prove the symmetrization identity (where \(\alpha+\beta=1\))
\begin{equation}\label{Symm} \begin{array}{l}\sum_{\sigma \in \mathbf{S}_{N}} \operatorname{sgn} \sigma\Big(\prod_{i<j}\left(\alpha+\beta \xi_{\sigma(i)} \xi_{\sigma(j)}-\xi_{\sigma(i)}\right)\\
\times \frac{\xi_{\sigma(2)} \xi_{\sigma(3)}^{2} \xi_{\sigma(4)}^{3} \cdots \xi_{\sigma(N)}^{N-1}}{\left(1-\xi_{\sigma(1)} \xi_{\sigma(2)} \xi_{\sigma(3)} \cdots \xi_{\sigma(N)}\right) \cdots\left(1-\xi_{\sigma(N-1)} \xi_{\sigma(N)}\right)\left(1-\xi_{\sigma(N)}\right)}\Big)  =\alpha^{N(N-1) / 2} \frac{\prod_{i<j}\left(\xi_{j}-\xi_{i}\right)}{\prod_{j}\left(1-\xi_{j}\right)}\end{array} \end{equation}

For later asymptotic analysis, the substitution $\xi_i = (z_i-1)^{-1}$ will  be more convenient.

\subsection{Matching of notation}
In this section, we match notation between various papers. Matching notation only requires looking at one lattice site. The figure below summarizes the various dualities. This paper uses only the \(q\)--Hahn and quantum \(q\)--Krawtchouk polynomials, but more information is provided for reference. 

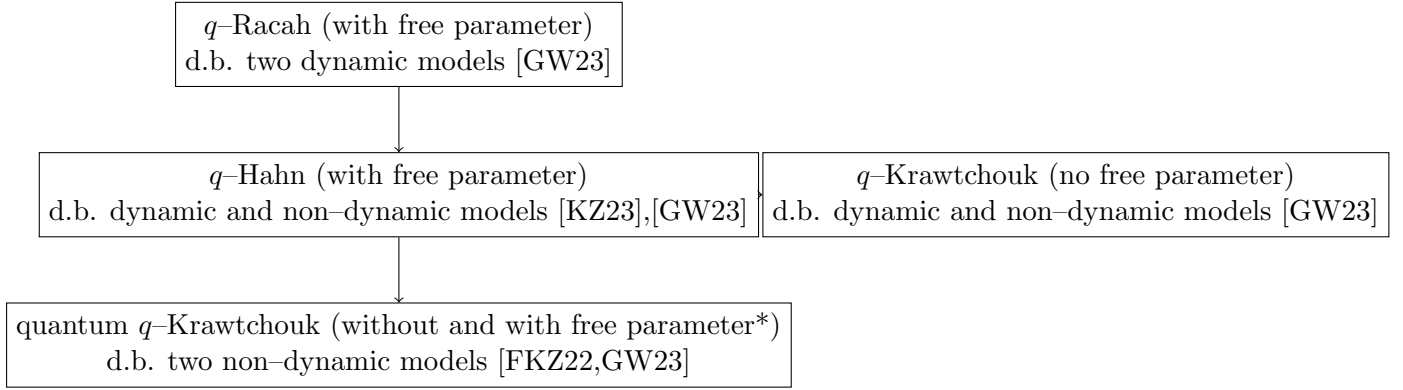
\begin{figure}[H]
\begin{tikzpicture}
\node at (4,6) [rectangle,draw,align=center] (T) {$q$--Racah (with free parameter) \\ d.b. two dynamic models [GW23]}; 
\node at (4,4) [rectangle,draw,align=center] (A) {$q$--Hahn (with free parameter) \\ d.b. dynamic and non--dynamic models [KZ23],[GW23]}; 
\node at (4,2) [rectangle,draw,align=center] (B) {quantum $q$--Krawtchouk (without and with free parameter*) \\ d.b. two non--dynamic models [FKZ22,GW23]}; 
\node at (13,4) [rectangle,draw,align=center] (C) {$q$--Krawtchouk (no free parameter)\\ d.b. dynamic and non--dynamic models [GW23]}; 
\draw [->] (T) -- (A);
\draw [->] (A) -- (B);
\draw [->] (A) -- (C);
\end{tikzpicture}
\caption{The various degenerations of named duality functions. The abbreviation ``d.b.'' means ``duality between.''}
\end{figure}



\subsubsection{Dualities from \cite{KZ23}}

The duality function is (where \(N\) is the number of lattice sites)

\normalsize
\begin{multline*}
\mathcal{D}_{c}(\mu, \xi)=\prod_{i=1}^{N}\left(e^{-2 \pi \mathrm{i} \lambda} q^{2\left(-c+N_{[1, i-1]}(\mu+\xi-2 J)+\mu_{i}-2 J_{i}\right)} ; q^{2}\right)_{2 J_{i}-\mu_{i}}\\
\times  { }_{3} \varphi_{2} \Big({q^{-2\left(2 J_{i}-\mu_{i}\right)}, q^{-2 \xi_{i}}, e^{-2 \pi \mathbf{i} \lambda} q^{2\left(2 N_{[1, i-1]}(\mu-J)+\mu_{i}-2 J_{i}\right)}}q^{-4 J_{i}}, e^{-2 \pi \mathbf{i} \lambda} q^{2\left(-c+N_{[1, i-1]}(\mu+\xi-2 J)+\mu_{i}-2 J_{i}\right)}; q^{2}, q^{2} \Big) \\
\times q^{-4 J_{i} N_{[1, i-1]}(\xi)-2 \xi_{i} N_{[1, i]}(\mu)}
\end{multline*}
\Large

If there is only one lattice site then we obtain 
\[
\left(e^{-2 \pi \mathrm{i} \lambda} q^{2\left(-c+\mu-2 J\right)} ; q^{2}\right)_{2 J-\mu}   { }_{3} \varphi_{2}\left(\begin{array}{c}q^{-2\left(2 J_{i}-\mu_{i}\right)}, q^{-2 \xi_{i}}, e^{-2 \pi \mathbf{i} \lambda} q^{2\left(\mu_{i}-2 J_{i}\right)}  \\q^{-4 J_{i}}, e^{-2 \pi \mathbf{i} \lambda} q^{2\left(-c+\mu_{i}-2 J_{i}\right)} \end{array} ; q^2, q^2\right)q^{-2 \xi_{i}\mu_i}
\]
which can be written in terms of the \(q\)--Hahn polynomial.

\subsubsection{Dualities from \cite{FKZ22}}
The duality function is given by the quantum \(q\)--Krawtchouk polynomial

\[
D_{\alpha_{i}}^{\boldsymbol{\theta}}\left(\boldsymbol{\xi}_{i}, \boldsymbol{\eta}_{i}\right)=\prod_{x=1}^{L} K^{\text{qtm}}_{\eta_{i}^{x}}\left(q^{-2 \xi_{i}^{x}}, p_{i}^{x}\left(\boldsymbol{\xi}_{i}, \boldsymbol{\eta}_{i}\right), \theta^{x}, q^{2}\right) 
\]
with
\[ p_{i}^{x}\left(\boldsymbol{\xi}_{i}, \boldsymbol{\eta}_{i}\right)=\alpha_{i}^{-1} q^{-2\left(N_{x-1}^{-}\left(\boldsymbol{\xi}_{i}\right)-N_{x+1}^{+}\left(\boldsymbol{\eta}_{i}\right)\right)+2 N_{x-1}^{-}(\boldsymbol{\theta})-1}.
\]
To match notation, we have that \(\theta = 2J.\) If there is only one lattice site, we have
\[
K_{\eta}^{\text{qtm}}(q^{-2\xi},\alpha^{-1}q^{-1},2J,q^2)
\]

\subsection{From duality to asymptotics}\label{fdta}

A previous proposition by the authors \cite{KZ23} gave a general method for proving asymptotics from orthogonal duality. The remainder of this subsection is essentially a verbatim repetition of that paper. 

Consider a filtered probability space \( (\Omega, \mathcal{F}, (\mathcal{F}_t)_{t \geq 0}, \mathbb{P}).\) Let \(X_t^{\lambda}\) be a stochastic process, depending on a parameter \(\lambda,\) with values in a countable state space \(\mathcal{S}.\) Let \(Q_{\lambda}(t)\) denote the probability measure on \(\mathcal{S}\) which is the pushforward of \(\mathbf{P}\) under \(X_t^{\lambda}.\) In other words, for any state \(x \in \mathcal{S},\) let \(Q_{\lambda}(t,x)\) denote the probability \( \mathbf{P}(X_t^{\lambda}=x).\)

Suppose that there is a family of measures, depending on \(\lambda\) on the state space \(\mathcal{S},\) which defines an inner product on the Hilbert space 
\[\mathcal{H}_{\lambda} = L^2(\mathcal{S},W_{\lambda}).\] In other words, 
\[
\langle f,g\rangle_{\lambda} = \sum_{x \in \mathcal{S}} f(x)g(x)W_{\lambda}(x).
\]
Suppose that the Hilbert space \(\mathcal{H}_{\lambda}\) has an orthogonal basis with respect to \(W_{\lambda},\) denoted by \( \{D_{\lambda}^S\}\) where \(S\) indexes the basis. In other words,
\begin{equation}\label{Ortho}
\sum_{x\in \mathcal{S}} D_{\lambda}^S(x)D_{\lambda}^{\overline{S}}(x)W_{\lambda}(x)=\delta_{S \overline{S}}\ n_{\lambda}^S
\end{equation}
for some non--negative normalization \(n_{\lambda}^S.\)

Let \(h(L)\) be a function on the state space \(\mathcal{S},\) and we assume that \(h(L)\) is an element of the Hilbert space \(\mathcal{H}_{\lambda}\) for all values of \(\lambda,\) but \(h(L)\) only depends on a large parameter \(L\) and not on \(\lambda.\) In the context of this paper, the letter \(h\) stands for ``height function.''

The previous paper states that if five inequalities hold for all values of \(\lambda,\) then the asymptotics of the height function  \(h(L)\) are the same for all values of \(\lambda.\) Below, we will assume that \(t=L\) for simplicity. If, for all \(S\) and \(L\), and all pairs \(\lambda,\overline{\lambda}\) we have
\begin{align}
\left(\sum_{y \in \mathcal{S}} h(L)[y]   \cdot \left|(n_{\lambda}^S)^{-1/2} D_{\lambda}^{S}(y) W_{\lambda}(y)-(n_{\overline{\lambda}}^S)^{-1/2} D_{\overline{\lambda}}^{S}(y) W_{\overline{\lambda}}(y)\right|\right) &\leq C_{1} L^{2} \label{oldeq1}\\
\sum_{y \in \mathcal{S}} h(L)[y] (n_{\lambda}^S)^{-1/2}D_{\lambda}^{S}(y) W_{\lambda}(y) &\leq C_{2} L^{2}\label{oldeq2}\\
\sum_{x \in \mathcal{S}}\left|(n_{\overline{\lambda}}^S)^{-1/2}D_{\overline{\lambda}}^{S}[x] Q_{\overline{\lambda}}(L, x)\right| &\leq M_{1}(L, S)\label{oldeq3}\\
\left(\sum_{x \in \mathcal{S}} (n_{\lambda}^S)^{-1/2}D_{\lambda}^{S}[x]\left(Q_{\lambda}(L, x)-Q_{\overline{\lambda}}(L, x)\right)\right) &\leq M_{2}(L, S)\label{oldeq4}
\end{align}
where
\[
\sum_{S} L^{2} M_{j}(L, S)<\infty, \quad \lim _{L \rightarrow \infty} L^{2} M_{j}(L, S)=0.
\]
Furthermore, assume

\begin{multline}\label{oldeq5}
\lim _{L \rightarrow \infty} \sum_{S}\left(\sum_{x \in \mathcal{S}}\left| (n_{\lambda}^S)^{-1/2}  D_{\lambda}^{S}[x]-(n_{\overline{\lambda}}^S)^{-1/2}D_{\overline{\lambda}}^{S}[x]\right| \cdot\left|Q_{\overline{\lambda}}(L, x)\right|\right)\\
\quad \quad \times \left(\sum_{y \in \mathcal{S}}\left|h(L)[y] (n_{\lambda}^S)^{-1/2}D_{\lambda}^{S}(y) W_{\lambda}(y)\right|\right)=0
\end{multline}
\Large
Then
\[
\lim _{L \rightarrow \infty} \sum_{x \in \mathcal{S}}\left|h(L)[x] Q_{\lambda}(L, x)-h(L)[x] Q_{\overline{\lambda}}(L, x)\right|=0
\]


\section{Main Results}

\subsection{From Orthogonal Duality to Determinants}

In this subsection, we significantly streamline the statements in section \ref{fdta}, at the cost of (slightly) losing generality.

\begin{proposition}\label{NewFDTA}
    Use the same notation as in section \ref{fdta}. Let \(\gamma(L)\) be some real--valued function. Assume that for all \(S\) and \(L\) and all pairs \(\lambda,\overline{\lambda}\) we have
    \begin{align}
         \sum_{y \in \mathcal{S}} \vert h(L)[y](n_{\lambda}^S)^{-1/2}D_{\lambda}^S(y) W_{\overline{\lambda}}(y) \vert & \leq  C\gamma(L) \label{Eq1}\\
         \sum_{x \in \mathcal{S}} \vert (n_{\lambda}^S)^{-1/2}D_{\lambda}^S(x)Q_{\overline{\lambda}}(L,x)\vert & \leq M(L,S) \label{Eq2}
    \end{align}
    where for each \(S,\) the function \( \gamma(L) M(L, S)\) is monotone in \(L\) and 
\(
\lim _{L \rightarrow \infty} \gamma(L) M(L, S)=0.
\)
Then
\[
\lim _{L \rightarrow \infty} \sum_{x \in \mathcal{S}}\left|h(L)[x] Q_{\lambda}(L, x)-h(L)[x] Q_{\overline{\lambda}}(L, x)\right|=0.
\]
\end{proposition}

\begin{remark}
The astute reader may notice that the normalization \(n_{\lambda}^S\) can be absorbed into the constants. Indeed, it is \textbf{orthogonality} of the duality functions that it signifcant, not the orthonormality. 
\end{remark}

\begin{remark}
In \eqref{Eq2}, if \(\lambda = \overline{\lambda}\) then the left hand side can be re--written as \(\mathbf{E}_x[D_{\lambda}^S(x(L))].\) By duality this equals \(\mathbf{E}_S[D_{\lambda}^{S(L)}(x)].\) This allows one to intuit its growth before any rigorous calculations.
\end{remark}

\subsection{Asympotics of dynamic ASEP}

The next theorem state asymptotics of dynamic ASEP, using the proposition in the previous subsection.

\begin{thm}\label{MainTheorem}
Let $s$ be any finite number and suppose that $\xi(t)$ being with step initial conditions at time $t=0$ (i.e. all lattice sites left of \(0\) are occupied, and all lattice site right of \(0\) are empty). Then, for $q>1$ and any \(\rho \in [-\infty,\infty)\), 
$$
\frac{h_s^+(\xi(t))}{F(t)} \rightarrow 0 
$$
in mean for all positive monotonic functions $F(t)$ such that $\lim_{t\rightarrow \infty } F(t) = \infty$.

When $\rho=\infty,$ one obtains the well--known Tracy--Widom \(F_2\) asymptotics of ASEP with step initial condition \cite{TW09,BCS12}.

\end{thm}

\subsection{Particle locations in (multi--species) ASEP}
This theorem is a generalization of a result in Theorem 2.1 of \cite{TW08}.

\begin{proposition}\label{TWProp}
Consider a $N$--particle ASEP with exactly one particle of each species. Let the initial condition consist of particles at $\mathbf{y}=(y_1<\ldots<y_N)$ with $q$--exchangeable distribution; in other words
\[
\mathrm{Prob}(\mathbf{y},\sigma) = \frac{q^{\inv(\sigma)}}{[N]_q^!}
\]
For $1 \leq i \leq N$, let $z_i(t)$ denote the location of the $i$th species particle at time $t$. For any fixed $K<N$ and any $x_1,\ldots,x_K$, let $\omega \in S(K)$ denote the permutation satisfying $x_{\omega(1)}< \ldots < x_{\omega(K)}$. Assume that $x_1,\ldots,x_K \geq M_{K+1} \geq \ldots \geq M_N$. Then
\begin{align}
&\mathrm{Prob}(z_1(t) = x_1, \ldots, z_K(t) = x_K, z_{K+1}(t) \geq M_{K+1} , \ldots, z_N(t) \geq M_N)\notag\\
&=  \frac{q^{\inv(\omega)}}{[K]^!_{q}} \sum_{\sigma \in S(N)} \left( \frac{1}{2\pi i}\right)^N \notag \\
&  \times \int_{C_r} \cdots \int_{C_r}   A_{\sigma} \prod_{i=1}^K \xi_{\sigma(i)}^{x_i - y_{\sigma(i)}-1} \prod_{i=K+1}^N \frac{ \xi_{\sigma(i)}^{M_i - y_{\sigma(i)}-1} }{ 1- \xi_{\sigma(i)}}  e^{\sum_{i=1}^N (p\xi_i^{-1}+q\xi_i-1)t} d\xi_1\cdots d\xi_N.
\end{align}
\end{proposition}
Note that setting $K=0$ and $M_1=\ldots=M_N$, we can symmetrize to recover Theorem 2.1 of \cite{TW08}, which gives the master equation for ASEP with \(N\) particles.

For the next proposition, we define the rainbow step initial conditions $\eta(0)=\{\eta_x(0):x \in \mathbf{Z}\}$  by 
$$
\eta_x(0) = 
\begin{cases}
x, \text{ if } x \geq 1 \\
0, \text{ if } x \leq 0.
\end{cases}
$$

\begin{proposition}\label{BlockDual}
Let \(M_1 \geq \cdots \geq M_n\) and consider some sequence of positive integers \(c_1 < \ldots < c_n.\) Let \((\mathbf{x}(0),\sigma(0))\) be rainbow step initial conditions. Let \(\mathbf{A}\) be the particle configuration with a particle of type \(c_j\) at \(M_j\) for all \(j.\) Then
\[
\mathbf{E}\left[\sum_{\tau \in S_n } \frac{q^{\inv(\tau)}}{[n]_q^!}
1_{ ((\mathbf{x}(t),\sigma(t)),\mathbf{A})\in I } \prod_{j=1}^n  q^{-2M_j + 2N^{c_{\tau{(j)}},-}_z( (\mathbf{x}(t),\sigma(t)))  }\right]
\]
equals the contour integral in Proposition \ref{TWProp} when \(K=0\). 
\end{proposition}

\begin{remark}
By using the color--blind Markov projection, we can consider the case when \(c_1=\cdots = c_n.\) This results a sum over \(\tau\in S_1,\) and then there are no permutations in the observable. 
\end{remark}

\section{Proofs}

\subsection{Proof of Proposition \ref{NewFDTA}}
This follows from proving that equations \eqref{Eq1} and \eqref{Eq2} imply equations \eqref{oldeq1}--\eqref{oldeq5}. Setting \(\lambda = \overline{\lambda}\) proves \eqref{oldeq2} and \eqref{oldeq3}, while taking distinct values proves \eqref{oldeq1} and \eqref{oldeq4}. Equation \eqref{oldeq5} then follows from \eqref{Eq1} , \eqref{Eq2} , and \eqref{oldeq1}--\eqref{oldeq4}. We replaced
\[
\sum_{S} \gamma(L)M(L,S)<\infty
\]
with the monotonicity condition to allow for an application of the monotone convergence theorem instead of the dominated convergence theorem.

\subsection{Proof of Propositions \ref{TWProp} and \ref{BlockDual}}

\subsubsection{Proof of Proposition \ref{BlockDual}}
The key observation here is that for a certain choice of \( \mathbf{x},\sigma,\) the duality function \(\widehat{D}\) becomes an indicator function on particle positions. 

Let \(\mathbf{A}(0)\) be the same as in the pposition. Note that for rainbow step initial conditions $\eta(0)$, the duality function $\widehat{D}(\eta(0),A(t))$ becomes an indicator function:
$$
\widehat{D}(\eta(0),\mathbf{A}(t)) = 
\begin{cases}
1 , &\text{ if } A^j(t) \subseteq [M_j,\infty) \text{ for all } j\\
0, &\text{ else }
\end{cases}
$$
Therefore, by the duality result,
$$
\mathbb{E}[ \widehat{D}(\eta(t),\mathbf{A}(0))  ] = \mathbb{E}[ \widehat{D}(\eta(0),\mathbf{A}(t))  ] = \mathbb{P}\left(    A^j(t) \subseteq [M_j,\infty)     \right).
$$
This completes the proof, pending the proof of Proposition \ref{TWProp} in the next section.


\subsubsection{Proof of Proposition \ref{TWProp} }
We first provide heuristics. We will take the summation of
\[
  q^{\inv(\tau)} \sum_{\sigma \in S_N} \int_{\mathcal{C}_r} \cdots \int_{\mathcal{C}_r} A_{\sigma} \prod_{i} \xi_{\sigma(i)}^{x_i - y_{\sigma(i)}-1} e^{(1+q)\sum_i \epsilon(\xi_i)t} d\xi_1 \cdots d\xi_N
\]
which we want to show equals
\[
  q^{\inv(\omega)}\sum_{\sigma \in S_N}   \int_{C_r} \cdots \int_{C_r}   A_{\sigma} \prod_{i=1}^K \xi_{\sigma(i)}^{x_i - y_{\sigma(i)}-1} \prod_{i=K+1}^N \frac{ \xi_{\sigma(i)}^{M_i - y_{\sigma(i)}-1} }{ 1- \xi_{\sigma(i)}}  e^{(1+q)\sum_i \epsilon(\xi_i)t} d\xi_1\cdots d\xi_N.
\]
If this summation is done naively, by ignoring the permutations and contour integrals, and allowing for two particles to occupy the same site, then the result would follow immediately from the geometric series
\[
\sum_{x_i \geq M_i} \xi^{x_i} = \frac{\xi^{M_i}}{1-\xi}.
\]
In a sense, the quantity \(\inv(\tau)\) ``accounts'' for exclusion in a way that this naive summation produces the correct answer.

Now note that for \(q=\beta/\alpha,\) where \(\alpha+\beta=1,\)
\begin{equation}\label{S1}
\frac{1}{1-\xi_i} - \frac{q\xi_j}{1-\xi_j} = \frac{1 - (1+q)\xi_j +q\xi_i\xi_j}{(1-\xi_i)(1-\xi_j)} =  \alpha^{-1}\frac{\alpha + \beta\xi_i \xi_j - \xi_j}{(1-\xi_i)(1-\xi_j)}.
\end{equation}
which matches one of the terms in the symmetrization identity \eqref{Symm}.
We first show that
\begin{multline}\label{Want}
\int_{C_r} \cdots \int_{C_r} \sum_{\sigma \in S(N)} A_{\sigma} \left[ \prod_{j=1}^{k-1} \xi_{\sigma(j)}^{\bar{x}_j - y_{\sigma(j)}-1} \right] \xi_{\sigma(k)}^{-M-y_{\sigma(k)}-1} \xi_{\sigma(k+1)}^{-M-y_{\sigma(k+1)}-1} \left[ \prod_{j=k+2}^N \xi_{\sigma(j)}^{\bar{x}_j - y_{\sigma(j)}-1} \right]   \\ 
\times \left( \frac{1}{1-\xi_{\sigma(k)}} - \frac{q\xi_{\sigma(k+1)}}{1-\xi_{\sigma(k+1)}}\right) e^{-\sum_i (1+q)\epsilon(\xi_i)t}  d\xi_1 \cdots d\xi_N =0. 
\end{multline}
We will see shortly that the \(q\) in the numerator arises from the quantity \(q^{\inv(\tau)}\) and that the contour integrals are not needed.

We partition $S(N)$ into two sets $S(N) = S^+_k(N) \cup S^-_k(N)$, where
$$
S^+_k(N) = \{\sigma \in S(N): \sigma(k)<\sigma(k+1)\}, \quad S^-_k(N) = \{\sigma \in S(N): \sigma(k)>\sigma(k+1)\}
$$
Noticing that $\sigma \in S_k^+(N)$ if and only if $\sigma^{(k)}:=\sigma \cdot (k \ k+1) \in S_k^-(N)$, and $\sigma(j) = \sigma^{(k)}(j)$ for all $j \neq k,k+1$, the integrand can be rewritten as (using \eqref{S1})
\begin{multline*}
\sum_{\sigma \in S^+_k(N)}  \left[ \prod_{j=1}^{k-1} \xi_{\sigma(j)}^{\bar{x}_j - y_{\sigma(j)}-1} \right] \frac{\xi_{\sigma(k)}^{-M-y_{\sigma(k)}-1} \xi_{\sigma(k+1)}^{-M-y_{\sigma(k+1)}-1} }{(1-\xi_{\sigma(k)})(1-\xi_{\sigma(k+1)})} \left[ \prod_{j=k+2}^N \xi_{\sigma(j)}^{\bar{x}_j - y_{\sigma(j)}-1} \right] \\
\times \alpha^{-1} \left(A_{\sigma}( \alpha + \beta\xi_{\sigma(k)}\xi_{\sigma(k+1)}- \xi_{\sigma(k+1)} ) + A_{\sigma^{(k)}}( \alpha + \beta\xi_{\sigma(k)}\xi_{\sigma(k+1)}- \xi_{\sigma(k)} )  \right) \\\times e^{-\sum_{j=1}^N (\alpha \xi_i^{-1}+\beta\xi_i-1)t}  
\end{multline*}
Noticing that $A_{\sigma^{(k)}} = A_{\sigma} S_{\sigma(k),\sigma(k+1)}$ for $\sigma \in S_k^+(N)$, we thus have 
$$
A_{\sigma}( \alpha + \beta\xi_{\sigma(k)}\xi_{\sigma(k+1)}- \xi_{\sigma(k+1)} ) + A_{\sigma^{(k)}}( \alpha + \beta\xi_{\sigma(k)}\xi_{\sigma(k+1)}- \xi_{\sigma(k)} ) = 0.
$$
This shows that \eqref{Want} is true.

Next, we use an additive property of \(\inv.\) Namely, let \( \tau \in S_N\) and let \(H\) be a normal subgroup of \(S_N.\) Then there is a unique decomposition \(\tau = a\omega\) such that \(\omega \in H\) and \(a\in S_N/H\) with \(\inv(\tau) = \inv(a) + \inv(\omega).\) We will apply this for \(H=S_K.\) 

We can now complete the proof. By \eqref{Want}, including values where \(x_k=x_{k+1}\) does not contribute to the summation. Therefore
\begin{align*}
&\mathrm{Prob}(z_1(t) = x_1, \ldots, z_K(t) = x_K, z_{K+1}(t) \geq M_{K+1} , \ldots, z_N(t) \geq M_N)\\
& = \sum_{j=K+1}^N\sum_{x_{j} \geq M_{j}}\sum_{\tau} \left( \frac{1}{2\pi i}\right)^N \frac{q^{\inv(\tau)}}{[N]_q^!}\\
& \quad \quad \times \sum_{\sigma \in S_N} \int_{\mathcal{C}_r} \cdots \int_{\mathcal{C}_r} A_{\sigma} \prod_{i\leq K}  \xi_{\sigma(i)}^{x_i - y_{\sigma(i)}-1} \prod_{i\geq K+1} \xi_{\sigma(i)}^{x_i - y_{\sigma(i)}-1} e^{(1+q)\sum_i \epsilon(\xi_i)t} d\xi_1 \cdots d\xi_N.
\end{align*}
The sum over \(\tau\) is actually a sum over \(a \in S_N/S_K\) where \(\tau=a\omega.\) Since
\[
\sum_{a \in S_N/S_K} q^{\inv(a)} = \frac{[N]_q^!}{[K]_q^!},
\]
using the geometric series completes the proof.

\subsection{Proof of Theorem \ref{MainTheorem}}
This is an application of Proposition \ref{NewFDTA}. For \(\rho=-\infty,\) the usual ASEP, we know that \(h_s^+(\xi(t))/F(t)\rightarrow 0\) almost surely \cite{BF87}, and Proposition \ref{NewFDTA} allows us to show convergence to \(0\) in mean. We need to show the two inequalities in the conditions of the proposition hold. First, we will match the notation, and then prove estimates on the duality function \(D,\) the weights \(W,w\) and the normalization \(n.\)

\subsubsection{Matching Notation}

Recall that in the notation \(D_{\lambda}^S(y),\) the symbol \(y\) denotes a state in the ``original'' process, while the symbol \(S\) denotes a state in the ``dual'' process. The parameter \(\lambda\) is the dynamic parameter, which we allow to take values in \([0,\infty].\) Recall that \(\rho\in (-\infty,\infty)\) is related to \(\lambda\) by \(\lambda = e^{\rho}.\)

Recall that the orthogonality relation can be written
\[
\sum_{\xi \in X} P_{\mathrm{R}}^{v}(\eta, \xi) P_{\mathrm{R}}^{v}\left(\eta^{\prime}, \xi\right) \omega_{\mathrm{R}}^{\mathrm{p}}(|\xi|) W_{\mathrm{R}}(\xi ; \vec{N}, \rho ; q)=\frac{\delta_{\eta, \eta^{\prime}}} {w(\eta ; \vec{N} ; q)    {\omega^{\mathrm{p}}(|\eta|)   }}
\]
and that there is the rescaling 
\[
K_{\mathrm{qtm}}^{v}(\eta, \xi)=\lim _{\rho \rightarrow\infty} \frac{\left(v^{-2} q^{-2 \rho}\right)^{|\eta|}}{c^{v}(|\eta|,|\xi| ; 0,0) C^{v}(|\eta|,|\xi| ; 0,2 \rho)} P_{\mathrm{R}}^{v q^{\rho}}(\eta, \xi).
\]
We then define
\[
D_{\lambda}^S(\xi) = \frac{\left(v^{-2} q^{-2 \rho}\right)^{|S|}}{c^{v}(|\eta|,|\xi| ; 0,0) C^{v}(|S|,|\xi| ; 0,2 \rho)} P_{\mathrm{R}}^{vq^{\rho}}(S,y) \]
and 
\[
W_{\lambda}(\xi) = \left(\frac{\left(v^{-2} q^{-2 \rho}\right)^{|S|}}{c^{v}(|\eta|,|\xi| ; 0,0) C^{v}(|S|,|\xi| ; 0,2 \rho)} \right)^{-2}W_{\mathrm{R}}(\xi ; \vec{N}, \rho ; q) \omega_{\mathrm{R}}^{\mathrm{p}}(\vert \xi \vert)\]
and
\[
n_{\lambda}(S)= \frac{1}{w(S ; \vec{N} ; q)    \omega^{\mathrm{p}}(|S|)   }
\]

With this notation, we have that
\[
\sum_{\xi \in X} D_{\lambda}^S(\xi)D_{\lambda}^{\overline{S}}(\xi)W_{\lambda}(\xi) = \delta_{S,\overline{S}}n_{\lambda}(S)
\]
which is equation \eqref{Ortho}, the orthogonality relationship necessary to apply Proposition \ref{NewFDTA}.

Also recall the explicit values
\[
\begin{aligned} \omega^{p}(x) & =\frac{v^{-2 x} q^{x(2 x-1)}\left(-v q^{\rho-|\vec{N}|+1} ; q^{2}\right)_{x}}{\left(v q^{-\rho+2 x-|\vec{N}|+1} ; q^{2}\right)_{|\vec{N}|-x}} \\ \omega_{\mathrm{R}}^{p}(x) & =\frac{\left(v q^{-\rho-2 x+|\vec{N}|+1} ; q^{2}\right)_{x}}{\left(-v q^{\rho-|\vec{N}|+1} ; q^{2}\right)_{x}}\end{aligned}
\]
and
\begin{multline*}
c^{v}(|\zeta|,|\xi| ; \lambda, \rho):=c(\zeta, \xi)\\
=\prod_{k=1}^{M} \frac{\left(v q^{2 \zeta_{k}-h_{k+1}^{+}(\xi)+h_{k-1}^{-}(\zeta)-N_{k}+1} ; q^{2}\right)_{N_{k}}}{\left(v q^{-2 \xi_{k}-h_{k+1}^{+}(\xi)+h_{k-1}^{-}(\zeta)+N_{k}+1} ; q^{2}\right)_{\xi_{k}+\zeta_{k}}}=\frac{\left(v q^{\lambda-\rho+2|\zeta|-|\vec{N}|+1} ; q^{2}\right)_{|\vec{N}|-|\zeta|}}{\left(v q^{\lambda-\rho-2|\xi|+|\vec{N}|+1} ; q^{2}\right)_{|\xi|}}
\end{multline*}
and
\[
C^{v}(|\zeta|,|\xi| ; \lambda, \rho)=\prod_{k=1}^{L} \frac{\left(-v q^{h_{k+1}^{+}(\xi)+h_{k-1}^{-}(\zeta)-N_{k}+1} ; q^{2}\right)_{\xi_{k}}}{\left(-v q^{h_{k+1}^{+}(\xi)+h_{k-1}^{-}(\zeta)-N_{k}+1} ; q^{2}\right)_{\zeta_{k}}}=\frac{\left(-v q^{\lambda+\rho-|\vec{N}|+1} ; q^{2}\right)_{|\xi|}}{\left(-v q^{\lambda+\rho-|\vec{N}|+1} ; q^{2}\right)_{|\zeta|}}
\]

For the proof, it is important the state \(S\) in the duality process only has finitely particles, so we write it in terms of particle variables
\[
S = \{s_1,\ldots,s_l\}
\]
where \(l\) is the number of particles in the particle configuration \(S.\) We write \(l=l(S)\) to emphasize the dependence on \(S.\) Then the duality function is expressed over occupied sites on \(S,\) or in other words: 
\[
P_{\mathrm{R}}^{v}(S, \xi)=\prod_{j=1}^{l(S)} p\left(1, \xi_{s_j} ; h_{s_j-1,0}^{-}(S), h_{s_j+1}^{+}(\xi), v, N_{k} ; q\right)
\]
here \(\xi\) is the original dynamic model and \(S\) is the dual model. 

\subsubsection{Rescaling as lattice size \(M\) diverges}
Note that the previous results of \cite{GW23} assumed a finite lattice size \(M,\) while in this paper we assume an infinite lattice size. The arguments for taking this limit are standard: one simply rescales the duality function by a constant. We briefly state that argument here. 

If the reader checks every term in the duality function, she will notice that the only dependence on the lattice size occur in the weights \(W_{\mathrm{R}}(\xi),\omega^{\mathrm{p}}(\vert S\vert),\omega_{\mathrm{R}}^{p}(\vert \xi \vert) \) and constants \(c^v,C^v.\)  The \(W_{\mathrm{R}}(\xi)\) term is more straightforward because \(h_{k+1}^+(\xi)\rightarrow \infty,\) so in the expressions
\[
W_{\mathrm{R}}(\xi ; q, \rho)  =\prod_{k\in \mathbb{Z}} W\left(\xi_k ; q, 1, h_{k+1}^{+}(\xi)\right),
\]
$$
\begin{aligned}
W(x ; q, 1, \rho) & =\frac{1+q^{4 x+2 \rho-2 }}{1+q^{2 \rho-2}} \frac{\left(-q^{2 \rho-2 } ; q^2\right)_x}{\left(-q^{2 \rho+2} ; q^2\right)_x} \frac{q^{-x(2 \rho+1+x-2 )}}{\left(-q^{-2 \rho} ; q^2\right)_1}\\
&=
\begin{cases}
    \frac{1}{1+q^{-2\rho}}, & x=0\\
    \left[\frac{1+q^{2\rho+2}}{1+q^{2\rho -2}} \frac{1+q^{2\rho -2}}{1+q^{2\rho +2}}\right] \frac{q^{-2\rho}}{1+q^{-2\rho}}, & x=1
\end{cases}
\end{aligned}
$$
the weight converges to a finite value with no rescaling in the \(q\)--Pochhammer in the denominator. (Note that in the \(x=1\) case the terms in the brackets cancel). The numerator requires a rescaling in \(q^M\) from the value of the height function \(h_{k+1}^+(\xi).\) Even more straightforwardly, 
\[
\omega_{\mathrm{R}}^{p}(\vert \xi \vert)  =\frac{\left(v q^{-\rho-2 \vert \xi \vert+|\vec{N}|+1} ; q^{2}\right)_{\vert \xi \vert}}{\left(-v q^{\rho-|\vec{N}|+1} ; q^{2}\right)_{\vert \xi \vert}} \rightarrow \textrm{const}
\]
by setting \( \vert \xi \vert = \vert \vec{N} \vert.\)

For the other terms, first note that
\begin{align*}
\omega^{\mathrm{p}}(\vert S\vert) &\sim \frac{v^{-2\vert S\vert}q^{\vert S\vert(2\vert S\vert-1)}}{(1-vq^{-\rho+2\vert S\vert-\vert \vec{N}\vert +1}) \cdots(1-vq^{-\rho+2\vert S\vert+2(\vert \vec{N}\vert - \vert S\vert-1) +1})   },\\
c^v(S, \xi,0,0) &\sim \frac{(1-vq^{ 2\vert S \vert - \vert\vec{N}\vert + 1 } ) \cdots ( 1-vq^{ \vert\vec{N}\vert + 1 }  ) }{(1-vq^{-2\vert \xi\vert + \vert \vec{N} \vert +1})\cdots ( 1-vq^{\vert \vec{N} \vert +1}   )}.
\end{align*}
where \( f \sim g\) means that \(\lim f/g\rightarrow 1.\) Note that the \(c^v\) term simplifies, because in the limit we can take \( \vert \xi \vert = \vert \vec{N} \vert.\) Then
\begin{multline*}
c^v(S, \xi,0,0) \sim \frac{1}{(1-vq^{-\vert \vec{N} \vert +1})\cdots (1-vq^{2 \vert {S}\vert-\vert \vec{N} \vert -1})   } \\
\sim (-v^{-1})^{ \vert S \vert+1}q^{(\vert S \vert+1)(\vert \vec{N}\vert)} q^{-(\vert S\vert+1)(\vert S\vert+2)/2 }
\end{multline*}
A similar calculation holds for \(\omega^{\mathrm{p}},\) with the key observation being that it does not depend on \(S\) except through \(\vert S\vert.\)

\subsubsection{Bounds on the weights \(\omega^{\mathrm{p}}_{\mathrm{R}}    ,W\)}
Recall that these are weights on \(\xi,\) the original dynamic process. As a reminder, 
\[
\omega_{\mathrm{R}}^{p}(\vert \xi \vert)  \leq \mathrm{const}.
\]
Furthermore, after re--normalization,
\begin{equation}\label{Winq}
W_{\mathrm{R}}(\xi;q,\rho) \leq \prod_{\substack{s \geq 0 \\ \xi_{s}=1}} q^{-2h_{s+1}^+(\xi)}
\end{equation}


\subsubsection{Bounds on the duality function}
Plugging in \(h_{s_j-1,0}^{-}(S)=j-1,\) we then have
\begin{align*}
P_{\mathrm{R}}^{v}(S, \xi)&=\prod_{j=1}^{l(S)} p\left(1, \xi_{s_j} ; j-1, h_{s_j+1}^{+}(\xi), v, 1 ; q\right)\\
&= \prod_{j=1}^{l(S)}  c_{\mathrm{p}}(1,\xi_{s_j}, j-1,h_{s_j+1}^{+}(\xi),v,1,q)\\
& \quad \quad \quad  \times P_{\xi_{s_j}}(q^{-2},-v q^{h_{s_j+1}^+(\xi)+j-3}, v^{-1} q^{h_{s_j+1}^+(\xi)-j+1},1;q^2)
\end{align*}
The constant 
\[
c_{\mathrm{p}}(n, x ; \lambda, \rho, v, N ; q)=v^{n} \frac{\left(-v q^{\rho+\lambda-N+1} ; q^{2}\right)_{x}\left(v q^{2 n-\rho+\lambda-N+1} ; q^{2}\right)_{N}}{q^{n(n+\rho+\lambda-N)}\left(v q^{-2 x-\rho+\lambda+N+1} ; q^{2}\right)_{x+n}}
\]
can be bounded above by 
\[
\mathrm{C}(v,S,\xi,\rho) q^{h_{s_j+1}^+(\xi)} 
\]
for some constant \(\mathrm{C}\) depending on the free parameter \(v\) and the states \(S,\xi\) and the dynamic parameter \(\rho.\)

Meanwhile, the \(q\)--Hahn polynomial can be expressed as
\begin{multline*}
P_{\xi_{s_j}}(q^{-2},-v q^{h_{s_j+1}^+(\xi)+j-3}, v^{-1} q^{h_{s_j+1}^+(\xi)-j+1},1;q^2) \\
= \ _{3}\varphi_2\left(q^{-2}, q^{-2\xi_{s_j}}, -q^{2h_{s_j+1}^+(\xi)+2j-6}q^{2h_{s_j+1}^+(\xi)-2j+2}q^{2\xi_{s_j}+2};-v q^{2h_{s_j+1}^+(\xi)+2j-6},q^{-2};q^2,q^2  \right)
\\
= \begin{cases}
    1 + \frac{(q^{-2};q^2)_1(q^{-2};q^2)_1\left(q^{4h_{s_j+1}^+(\xi)};q^2\right)_1}{(q^{-2};q^2)_1\left(-vq^{2h_{s_{j}+1}^+(\xi)+2j-6};q^2\right)_1}\frac{q^2}{(q^2;q^2)_1} & \quad \xi_{s_j}=1, \\
    1, & \quad \xi_{s_j}=0. 
\end{cases}
\end{multline*}
Therefore, define the function, while recalling that the height function \(h\) depends on \(\rho,\)
\[
B(s,j,v,\xi,\rho) = 1 - \frac{1-q^{4h_{s+1}^+(\xi)}}{1+vq^{2h^+_{s+1}(\xi)+2j-6}}.
\]
We obtain the bound
\[
D_{\lambda}^S(\xi)  \leq \prod_{\substack{1 \leq j \leq l(S) \\ \xi_{s_j}=1 }} \mathrm{C}(v,\vert S\vert,\vert \xi\vert,\rho) q^{h_{s_j+1}^+(\xi)} B(s,j,v,\xi,\rho)
\]
where \(\mathrm{C}(v,\vert S\vert,\vert \xi\vert,\rho)\) is some constant. So finally
\begin{equation}\label{BCSIneq}
\left| D_{\lambda}^S(\xi)  \right| \leq \mathrm{\overline{C}}(v,\vert S\vert,\vert \xi\vert,\rho) \left[\prod_{\substack{1 \leq j \leq l(S) \\ \xi_{s_j}=1}} q^{2h_{s_j+1}^+(\xi)}\right]
\end{equation}
for some constant \( \overline{C}.\)


\subsubsection{Proof of \eqref{Eq1}}
Set \(\gamma(L)=1/F(L).\) We will prove
\[
 \sum_{y \in \mathcal{S}} \vert(n_{\lambda}^S)^{-1/2} h(L)[y]D_{\lambda}^S(y) W_{\overline{\lambda}}(y) \vert  \leq  C\gamma(L).
\]
Plugging in the above bounds \eqref{BCSIneq} and \eqref{Winq} for the duality function and weights, respectively, and absorbing \( (n_{\lambda}^S)^{-1/2}\) into the constant \(C\) we see that it suffices to show
\[
\sum_{\xi} \frac{h_r^+(\xi)}{F(L)} \left[\prod_{\substack{1 \leq j \leq l(S) \\ \xi_{s_j}=1}} q^{2h_{s_j+1}^+(\xi)}\right]  \left[\prod_{\substack{s \geq 0 \\ \xi_{s}=1}} q^{-2h_{s+1}^+(\xi)}\right] \leq \frac{A}{F(L)}
\]


\subsubsection{Proof of \eqref{Eq2}}
We want to prove
\[
         \sum_{x \in \mathcal{S}} \vert (n_{\lambda}^S)^{-1/2}D_{\lambda}^S(x)Q_{\overline{\lambda}}(L,x)\vert  \leq M(L,S) 
\]
where \(M(L, S)/F(L)\) is monotone in \(L\) and 
\[
 \lim _{L \rightarrow \infty}  M(L, S)/F(L)=0.
\]
We will show that taking \(M(L,S)\) to be a sufficiently large constant will work.

The inequality in \eqref{BCSIneq} gives an upper bound on the duality function \(D_{\lambda}^S(x);\) in fact the upper bound is the duality function from \cite{BCS12}, which we denote \( D_{\text{BCS}}(S,\xi).\) This is a duality function for the usual ASEP and its space reversal. Since $Q_{\overline{\lambda}}(L,\xi)$ is the weight for \(\xi\) at time \(L\) in the dynamic ASEP, it thus suffices to bound
\[
\mathbf{E}_{\text{step}}[D_{\text{BCS}}(S,\xi(L))].
\]
Because the dynamic parameter has the effect of pushing the height function down, we can couple the dynamic ASEP \(\xi(L)\) with the usual ASEP \(\xi'(L)\)
such that 
\[
\mathbf{E}_{\text{step}}[D_{\text{BCS}}(S,\xi(L))] \leq \mathbf{E}_{\text{step}}[D_{\text{BCS}}(S,\xi'(L))].
\]
By the duality relation, the term on the right--hand--side equals
\[
\mathbf{E}_{S}[D_{\text{BCS}}(S(L),\xi'(0))]
\]
where \(S(L)\) is an ASEP starting at initial condition \(S\) with drift to the right. Since \(\xi'(0)\) is step initial conditions with particles to the left of the origin, we thus obtain a constant as an upper bound. This completes the proof.

\section{Computer Simulations}

An animation of dynamic ASEP can be found in the first author's IPAM talk; slides are downloadable from \cite{TSB2}.

\begin{figure}[H]

\begin{center}
\includegraphics[height=5cm]{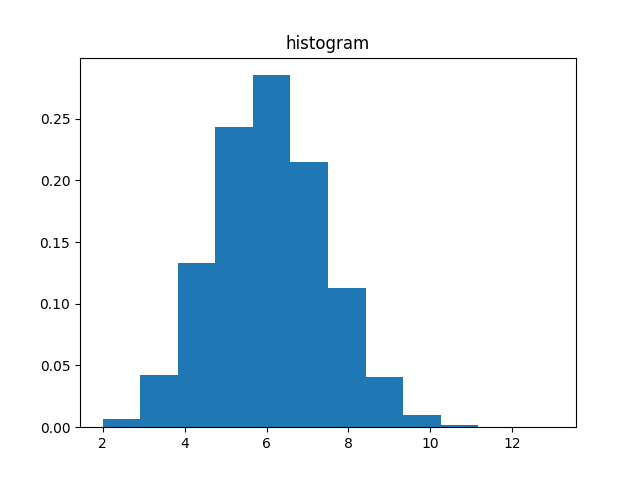} \includegraphics[height=5cm]{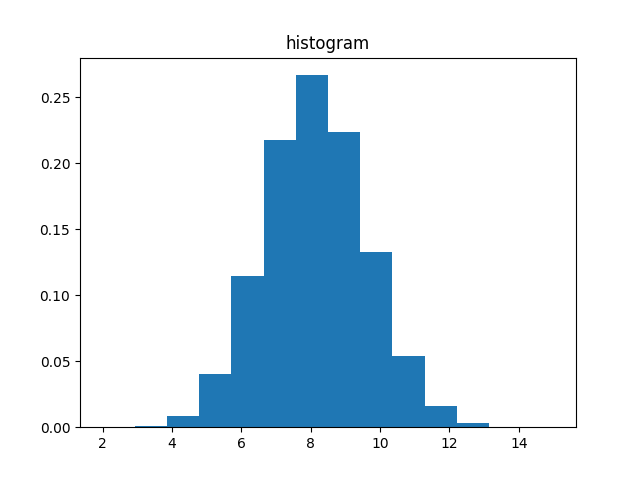}
\end{center}

\caption{Both figures on the top show 5000 samples of the height function \(h_x(t)\) of dynamic ASEP at time \(t=1000\) and spatial parameter \(x=0\). For both figures, the value of the dynamic parameter is \(1.1,\) while the asymmetry parameter \(q\) is \(0.9\) for the figure on the left and \(1.1\) for the figure on the right.}

\end{figure}

\begin{figure}[H]

\begin{center}
\includegraphics[height=5cm]{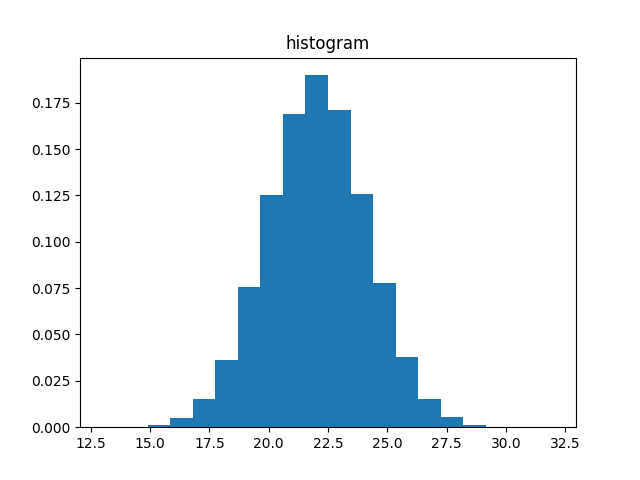}
\end{center}

\caption{In contrast, the ASEP for \(q=0.9\) has fluctuations which have larger values.}

\end{figure}

\end{document}